\documentclass[12pt]{article}

\usepackage{amssymb,amsmath,amsthm}

\newcommand{\ovl}{\overline}

\newcommand{\vp}{\varepsilon}
\newcommand{\cl}[1]{{\mathcal{#1}}}
\newcommand{\bb}[1]{{\mathbb{#1}}}
\newcommand{\tr}{\text{tr}}

\numberwithin{equation}{section}

\theoremstyle{plain}
\newtheorem{lem}{Lemma}[section]
\newtheorem{pro}[lem]{Proposition}
\newtheorem{thm}[lem]{Theorem}
\newtheorem{cor}[lem]{Corollary}

\theoremstyle{definition}

\theoremstyle{remark}
\newtheorem{rem}[lem]{Remark}

\begin{document}

\begin{center}
\LARGE UNITARY PERTURBATIONS OF MASAS\\
IN TYPE $II_1$ FACTORS
\end{center}\bigskip

{\hspace{1in}
\begin{tabular}{cc}
Allan M. Sinclair&Roger R.~Smith$^*$ \\
\noalign{\medskip}
Department of Mathematics&Department of Mathematics\\
University of Edinburgh&Texas A\&M University\\
Edinburgh, EH9 3JZ&College Station, TX \ 77843\\
SCOTLAND&U.S.A.\\
e-mail:\ {\tt allan@maths.ed.ac.uk}&{\tt rsmith@math.tamu.edu}
\end{tabular}
}
\vspace{.5in}

\abstract{
The main result of this paper is the inequality
\[
d(u,N(\cl A))/31 \le \|\bb E_{\cl A} - \bb E_{u\cl A u^*}\|_{\infty,2} \le
4d(u,N(\cl A)),
\]
where $\cl A$ is a masa in a separably acting type $II_1$ factor $\cl N$,
$u\in\cl N$ is a unitary, $N(\cl A)$ is the group of normalizing unitaries,
$d$ is the distance measured in the $\|\cdot\|_2$-norm, and
$\|\cdot\|_{\infty,2}$ is a norm defined on the space of bounded maps on $\cl
N$ by
\[
\|\phi\|_{\infty,2} = \sup\{\|\phi(x)\|_2\colon \ \|x\|\le 1\}.
\]
This result implies that a unitary which almost normalizes a masa must be close
to a normalizing unitary. The inequality also shows that every singular masa
is $\alpha$-strongly singular for $\alpha = 1/31$.}

\vfill

$\underline{\hspace{1in}}$

\noindent $^*$Partially supported by a grant from the National Science
Foundation.
\newpage

\section{Introduction}\label{sec1}

\indent

This paper is a continuation of our work in \cite{SS1,SS2}, where we
introduced and studied strong singularity for maximal abelian self-adjoint
subalgebras (masas) in a type $II_1$ factor $\cl N$. If $\phi\colon \ \cl N\to
\cl N$ is a bounded linear map, then $\|\phi\|_{\infty,2}$ denotes the quantity
\begin{equation}\label{eq1.1}
\|\phi\|_{\infty,2} = \sup\{\|\phi(x)\|_2\colon \ \|x\|\le 1\},
\end{equation}
and $\|\bb E_{\cl A} - \bb E_{\cl B}\|_{\infty,2}$ measures the distance
between two masas $\cl A$ and $\cl B$, where $\bb E_{\cl A}$ and $\bb E_{\cl
B}$ are the associated trace preserving 
conditional expectations. The singular masas are those
which contain their groups of unitary normalizers, \cite{Di}. Within this
class, we introduced the notion of strongly singular masas, \cite{SS1}. The
defining inequality for strong singularity of a masa $\cl A\subseteq \cl N$ is
\begin{equation}\label{eq1.2}
\|u - \bb E_{\cl A}(u)\|_2 \le \|\bb E_{\cl A} - \bb E_{u\cl A
u^*}\|_{\infty,2}
\end{equation}
for all unitaries $u\in \cl N$. Such an inequality implies that each
normalizing unitary lies in $\cl A$, so (\ref{eq1.2}) can only hold for
singular masas. We may weaken (\ref{eq1.2}) by inserting a constant $\alpha
\in (0,1]$ on the left hand side, and a masa which satisfies the modified
inequality is called $\alpha$-strongly singular. All the examples of singular
masas in \cite{SS1,SS2} are strongly singular, and so it is natural to ask
whether all singular masas have this property, or whether all are
$\alpha$-strongly singular (with $\alpha>0$ perhaps depending on the
particular singular masa). The first question remains open, but we will give a
positive answer to the second for the value $\alpha = 1/31$. This will be a
consequence of a much stronger inequality, valid for all masas, which is the
main result of the paper (Theorem \ref{thm5.3}). Due to certain technical
considerations, we will restrict attention throughout to type $II_1$ factors
which have separable preduals or, equivalently, have representations on
separable Hilbert spaces. If $u\in\cl N$ is unitary and $v\in N(\cl A)$, then
the inequality
\begin{equation}\label{eq1.3}
\|\bb E_{\cl A} - \bb E_{u\cl A u^*}\|_{\infty,2} \le 4\|u-v\|_2
\end{equation}
is straightforward (see Lemma \ref{thm5.2}). The converse, proved in
Theorem~\ref{thm5.1}, gives the existence of a unitary $v\in N(\cl A)$
satisfying
\begin{equation}\label{eq1.4}
\|u-v\|_2 \le 31\|\bb E_{\cl A} - \bb E_{u\cl Au^*}\|_{\infty,2}.
\end{equation}
This is a much deeper inequality, and most of the paper is devoted to its
proof. It is surely not the case that 31 is the best possible constant in
(\ref{eq1.4}), but this number emerges from a chain of various estimates and
we have been unable to improve upon it.

Strong singularity of a masa $\cl A$ was known, \cite{SS1}, to be related 
to the invariant $\delta(\cl A)$ which had been introduced by Popa in 1983, 
\cite{P0}. Recently, \cite{P3}, he showed that $\delta(\cl A)=1$ for all singular
masas $\cl A$, and so $\delta$ only takes the values $0$ and $1$. This result,
\cite[Cor. 2]{P3}, may be stated as follows. If $\cl A$ is a singular masa  in a type
$II_1$ factor $\cl N$ and $v$ is a partial isometry in $\cl N$ with 
$vv^*$ and $v^*v$ in $\cl A$, then
\begin{equation}\label{eq1.5}
\|vv^*\|_2={\mathrm{sup}}\,\{\|x-{\bb{E}}_{\cl A}(x)\|_2:\ x\in v{\cl A}v^*,\ \|x\|
\leq 1\}.
\end{equation}
This result supports the possibility that all singular masas are strongly singular. 
The method of proof in \cite{P3} (and of the main technical lemma in \cite{P4}) uses 
the method developed by Christensen, \cite{Ch}, together with the pull--down identity 
of $\Phi$ from \cite{PP} and some work by Kadison on center--valued traces, \cite{K}. 
Our main proof (Theorem \ref{thm4.2}) applies Christensen's technique to 
$\langle {\cl N},e_{\cl B}\rangle$, uses the pull--down map $\Phi$, and also requires 
approximation of finite projections in $L^{\infty}[0,1]\overline{\otimes}B(H)$. These 
are combined with a detailed handling of various inequalities involving 
projections and partial isometries.

There are two simple ways in which masas $\cl A$ and $\cl B$ in a type $II_1$ factor 
can be close in the $\|\cdot\|_{\infty,2}$--norm on their conditional expectations. If $u$
is a unitary close to $\cl A$ in $\|\cdot\|_2$--norm and $u{\cl A}u^*={\cl B}$, then
$\|\bb E_{\cl A}-\bb E_{\cl B}\|_{\infty,2}$ is small. 
Secondly, if there is a projection $q$ 
of large trace in $\cl A$ and $\cl B$ with $q{\cl A}=q{\cl B}$, then again
$\|\bb E_{\cl A}-\bb E_{\cl B}\|_{\infty,2}$ is small. In Theorem \ref{thm5.4} we show 
that a combination of these two methods is the only way in which $\cl A$ can be close to 
$\cl B$ in separably acting factors. In the last section, we will relate the distance between masas
to inclusions ${\cl A} \subset_{\delta} {\cl B}$ or 
${\cl B} \subset_{\delta} {\cl A}$ introduced in \cite{MvN} and studied in \cite{Ch}.
This does not hold for general von~Neumann subalgebras of $\cl N$ and is a special property 
of masas.

The crucial estimates are contained in Theorem \ref{thm4.2}, Proposition
\ref{pro3.4} and Corollary \ref{cor2.5}; all other results in those sections
are really parts of their proofs, broken down to manageable size. We recommend
reading these three results in the order stated, referring back to ancillary
lemmas and propositions as needed. Corollary~\ref{cor2.5} is essentially due
to Christensen, \cite{Ch}, but without the norm inequalities which we have
included. While we do not quote directly from 
\cite{P3,P4}, we note that the calculations in these papers have served to
motivate strongly the present work and to guide us through considerable
complications. We thank Sorin Popa for providing us with copies of these 
two recent preprints prior to their publication.
\newpage

\section{Preliminaries}\label{sec2}

\indent

Let $\cl N$ be a fixed but arbitrary separably acting type $II_1$ factor with
faithful normalized normal trace tr, and let $\cl B$ be a maximal abelian
subalgebra (masa) of $\cl N$. The trace induces an inner product
\begin{equation}\label{eq2.1}
\langle x,y\rangle = \tr(y^*x),\qquad x,y\in\cl N,
\end{equation}
on $\cl N$. Then $L^2(\cl N,\tr)$ is the resulting completion with norm
\begin{equation}\label{eq2.2}
\|x\|_2 = (\tr(x^*x))^{1/2},\qquad x\in\cl N,
\end{equation}
and when $x\in\cl N$ is viewed as a vector in this Hilbert space we will
denote it by $\hat x$. The unique trace preserving conditional expectation
$\bb E_{\cl B}$ of $\cl N$ onto $\cl B$ may be regarded as a projection in
$B(L^2(\cl N), \tr)$, where we denote it by $e_{\cl B}$. Thus
\begin{equation}\label{eq2.3}
e_{\cl B}(\hat x) = \widehat{\bb E_{\cl B}(x)},\qquad x\in\cl N.
\end{equation}
Properties of the trace show that there is a conjugate linear isometry
$J\colon \ L^2(\cl N,\tr)\to L^2(\cl N,\tr)$ defined by
\begin{equation}\label{eq2.4}
J(\hat x) = \widehat{x^*},\qquad x\in\cl N,
\end{equation}
and it is standard that $\cl N$, viewed as an algebra of left multiplication
operators on $L^2(\cl N,\tr)$, has commutant $J\cl NJ$. The von~Neumann algebra
generated by $\cl N$ and $e_{\cl B}$ is denoted by $\langle \cl N,e_{\cl
B}\rangle$, and has commutant $J\cl BJ$. Thus $\langle\cl N,e_{\cl B}\rangle$
is a type $I_\infty$ von~Neumann algebra, since its commutant is abelian, and
its center is $J\cl BJ$, a masa in $\cl N\,'$ and thus isomorphic to
$L^\infty[0,1]$. The general theory of type $I$ von~Neumann algebras,
\cite{KR}, shows that there is a separable infinite dimensional Hilbert space
$H$ so that $\langle\cl N,e_{\cl B}\rangle$ and $L^\infty[0,1]\overline\otimes
B(H)$ are isomorphic. When appropriate, we will regard an element $x\in\langle
\cl N,e_{\cl B}\rangle$ as a uniformly bounded measurable $B(H)$-valued
function $x(t)$ on [0,1]. Under this identification, the center $J\cl BJ$ of
$\langle\cl N,e_{\cl B}\rangle$ corresponds to those functions taking values
in $\bb CI$. We denote by $\tau$ the unique semi-finite faithful normal trace
on $B(H)$ which assigns the value 1 to each rank 1 projection.

The following lemma (see \cite{PP,P2}) summarizes some of the basic properties
of $\langle \cl N, e_{\cl B}\rangle$ and $e_{\cl B}$.

\begin{lem}\label{lem2.1}
Let $\cl N$ be a separably acting type $II_1$ factor with a masa $\cl B$. Then
\begin{align}
\label{eq2.5}
\text{\rm (i)}\quad &e_{\cl B} x e_{\cl B} = e_{\cl B}\bb E_{\cl B}(x) = \bb
E_{\cl B}(x) e_{\cl B},\qquad x\in\cl N;\\
\label{eq2.6}
\text{\rm (ii)}\quad &e_{\cl B}\langle \cl N,e_{\cl B}\rangle = e_{\cl B}\cl
N,\quad \langle \cl N,e_{\cl B}\rangle e_{\cl B} = \cl Ne_{\cl B};\\
\label{eq2.7}
\text{\rm (iii)}\quad &\text{if } x\in \langle\cl N,e_{\cl B}\rangle \text{
and } e_{\cl B}x = 0, \text{ then } x=0;\\
\label{eq2.8}
\text{\rm (iv)}\quad &e_{\cl B}\langle \cl N,e_{\cl B}\rangle e_{\cl B} = \cl
Be_{\cl B} = e_{\cl B}\cl B;\\
\text{\rm (v)}\quad &\text{there is a faithful normal semi-finite trace {\rm
Tr}  on $\langle\cl N,e_{\cl B}\rangle$ which satisfies}\nonumber
\end{align}
\begin{equation}\label{eq2.9}
\text{\rm Tr}(xe_{\cl B}y) = \text{\rm tr}(xy),\qquad x,y\in\cl N,
\end{equation}
and in particular,
\begin{equation}\label{eq2.10}
\text{\rm Tr}(e_{\cl B}) = 1.
\end{equation}
\end{lem}

The following result will be needed subsequently. We denote by
$\|x\|_{\text{Tr},2}$ the Hilbert space norm induced by Tr on the subspace of
$\langle\cl N, e_{\cl B}\rangle$ consisting of elements satisfying
$\text{Tr}(x^*x) < \infty$.

\begin{lem}\label{lem2.2}
Let $\vp>0$. If $f\in \langle \cl N,e_{\cl B}\rangle$ is a projection of
finite trace and
\begin{equation}\label{eq2.11}
\|f-e_{\cl B}\|_{\text{\rm Tr},2} \le \vp,
\end{equation}
then there exists a central projection $z\in \langle\cl N, e_{\cl B}\rangle$
such that $zf$ and $ze_{\cl B}$ are equivalent projections in $\langle\cl N,
e_{\cl B}\rangle$. Moreover, the following inequalities hold:
\begin{equation}\label{eq2.12}
\|zf-ze_{\cl B}\|_{\text{\rm Tr},2}, \quad \|ze_{\cl B}-e_{\cl
B}\|_{\text{\rm Tr},2},\quad \|zf-e_{\cl B}\|_{\text{\rm Tr},2} \le \vp.
\end{equation}
\end{lem}

\begin{proof}
From (\ref{eq2.8}), $e_{\cl B}$ is an abelian projection in $\langle\cl N,
e_{\cl B}\rangle$ so, altering $e_{\cl B}(t)$ on a null set if necessary, each
$e_{\cl B}(t)$ is a projection in $B(H)$ whose rank is at most 1. If $\{t\in
[0,1]\colon \ e_{\cl B}(t) = 0\}$ were not a null set, then there would exist
a non-zero central projection $p$ corresponding to this set so that $e_{\cl
B}p = 0$, contradicting (\ref{eq2.7}). Thus we may assume that each $e_{\cl
B}(t)$ has rank 1.

Since Tr  is a faithful normal semi-finite trace on $\langle\cl N, e_{\cl
B}\rangle$, there exists a non-negative $\bb R$-valued measurable function
$k(t)$ on [0,1] such that
\begin{equation}\label{eq2.13}
\text{Tr}(y) = \int^1_0 k(t) \tau(y(t))dt, \quad y\in \langle \cl N,e_{\cl
B}\rangle, \quad \text{Tr}(y^*y) < \infty.
\end{equation}
By (\ref{eq2.10}),
\begin{equation}\label{eq2.14}
\text{Tr}(e_{\cl B}) = \int^1_0 k(t) \tau(e_{\cl B}(t)) dt = 1,
\end{equation}
and thus integration against $k(t)$ defines a probability measure $\mu$ on
[0,1] such that
\begin{equation}\label{eq2.15}
\text{Tr}(y) = \int^1_0 \tau(y(t)) d\mu(t),\quad y\in \langle \cl N, e_{\cl
B}\rangle, \quad \text{Tr}(y^*y) < \infty.
\end{equation}
It follows from (\ref{eq2.15}) that
\begin{equation}\label{eq2.16}
\|y\|^2_{\text{Tr},2} = \int^1_0 \tau(y(t)^* y(t)) d\mu(t),\qquad y\in
\langle\cl N, e_{\cl B}\rangle.
\end{equation}

Consider a rank 1 projection $p\in B(H)$ and a projection $q\in B(H)$ of rank
$n\ge 2$. Then 
\begin{align}
\tau((p-q)^2) &= \tau(p+q-2pq)\nonumber\\
&= \tau(p+q-2pqp)\nonumber\\
&\ge \tau(p+q-2p)\nonumber\\
\label{eq2.17}
&\ge 1,
\end{align}
and the same inequality is obvious if $q=0$. Let $G = \{t\in [0,1]\colon \
\text{rank}(f(t))\ne 1\}$. Then, from (\ref{eq2.11}),
\begin{align}
\vp^2 &\ge \|f-e_{\cl B}\|^2_{\text{Tr},2}\nonumber\\
&\ge \int_G \tau((f(t)-e_{\cl B}(t))^2) d\mu(t)\nonumber\\
\label{eq2.18}
&\ge \mu(G),
\end{align}
by (\ref{eq2.17}). Let $z = \chi_{G^c} \otimes I$, a central projection in
$\langle\cl N, e_{\cl B}\rangle$. Then the ranks of $z(t)f(t)$ and
$z(t)e_{\cl B}(t)$ are simultaneously 0 or 1, and so $zf$ and $ze_{\cl B}$ are
equivalent projections in $\langle\cl N,e_{\cl B}\rangle$. Then
\begin{align}
\|ze_{\cl B}-e_{\cl B}\|^2_{\text{Tr},2} &= \int_G \tau(e_{\cl B}(t))
d\mu(t)\nonumber\\
&= \mu(G)\nonumber\\
\label{eq2.19}
&\le \vp^2,
\end{align}
from (\ref{eq2.18}), while
\begin{align}
\|zf-e_{\cl B}\|^2_{\text{Tr},2} &= \int_G \tau(e_{\cl B}(t))d\mu(t) +
\int_{G^c} \tau((f(t) - e_{\cl B}(t))^2) d\mu(t)\nonumber\\
\label{eq2.20}
&\le \|f-e_{\cl B}\|^2_{\text{Tr},2}
\end{align}
since, on $G$,
\begin{equation}\label{eq2.21}
\tau(e_{\cl B}(t)) = 1 \le \tau((f(t) - e_{\cl B}(t))^2),
\end{equation}
by (\ref{eq2.17}). Finally, 
\begin{equation}\label{eq2.22}
\|zf-ze_{\cl B}\|_{\text{Tr},2} \le \|z\| \|f-e_{\cl B}\|_{\text{Tr},2} \le
\vp, 
\end{equation}
completely the proof of (\ref{eq2.12}).
\end{proof}

We now recall some properties of the polar decomposition and some trace norm
inequalities. These may be found in \cite{Co,KR}.

\begin{lem}\label{lem2.3}
Let $\cl M$ be a von Neumann algebra.
\begin{itemize}
\item[\rm (i)] If $w\in\cl M$ then there exists a partial isometry  $v\in\cl
M$, whose initial and final spaces are respectively the closures of the ranges
of $w^*$ and $w$, satisfying
\begin{equation}\label{eq2.23}
w = v(w^*w)^{1/2} = (ww^*)^{1/2}v.
\end{equation}
\item[\rm (ii)] Suppose that $\cl M$ has a faithful normal semifinite trace
\text{\rm Tr}. If $x\in \cl M$, $0\le x\le 1$, \text{\rm Tr}$(x^*x) < \infty$,
and $f$ is the spectral projection of $x$ corresponding to the interval
{\rm [1/2, 1]}, then 
\begin{equation}\label{eq2.24}
\|e-f\|_{\text{\rm Tr},2} \le 2\|e-x\|_{\text{\rm Tr},2}
\end{equation}
for any projection $e\in\cl M$ of finite trace.
\item[\rm (iii)]
Suppose that $\cl M$ has a faithful normal semifinite trace and let $p$ and
$q$ be equivalent finite projections in $\cl M$. Then there exists a partial
isometry $v\in \cl M$ and a unitary $u\in\cl M$ satisfying
\begin{gather}\label{eq2.25}
v^*v = p,\quad vv^* = q,\\
\label{eq2.26}
v|p-q|  = |p-q|v,\\
\label{eq2.27}
|v-p|, |v-q| \le 2^{1/2} |p-q|,\\
\label{eq2.28}
upu^* = q,\quad u|p-q| = |p-q|u,\\
\label{eq2.29}
|1-u| \le 2^{1/2} |p-q|.
\end{gather}
\end{itemize}
\end{lem}
\medskip

The following result is essentially in \cite{Ch}, and is also used in 
\cite{P3,P4}. We reprove it here since the norm estimates that we obtain will
be crucial for subsequent developments.

\begin{pro}\label{pro2.4}
Let $\cl A$ and $\cl B$ be masas in a separably acting type $II_1$ factor $\cl
N$, and let $\ovl K^w_{\cl A}(e_{\cl B})$ be the weak closure of the set
\begin{equation}\label{eq2.30}
K_{\cl A}(e_{\cl B}) = 
{\mathrm{conv}}\,\{ue_{\cl B}u^*\colon \ u \text{ is a unitary in } \cl
A\}
\end{equation}
in $\langle\cl N,e_{\cl B}\rangle$. Then $\ovl K^w_{\cl A}(e_{\cl B})$
contains a unique element $h$ of minimal $\|\cdot\|_{\text{\rm Tr},2}$-norm,
and this element satisfies
\begin{align}
\label{eq2.31}
\text{\rm (i)}\quad &h\in \cl A'\cap \langle\cl N, e_{\cl B}\rangle,\qquad 0
\le h \le 1;\hspace{2.56in}\\
\label{eq2.32}
\text{\rm (ii)}\quad &1-\text{\rm Tr}(e_{\cl B}h) \le \|(I-\bb E_{\cl B})\bb
E_{\cl A}\|^2_{\infty,2};\\
\label{eq2.33}
\text{\rm (iii)}\quad &\text{\rm Tr}(e_{\cl B}h) = \text{Tr}(h^2);\\
\label{eq2.34}
\text{\rm (iv)}\quad &\|h-e_{\cl B}\|_{\text{\rm Tr},2} \le \|(I-\bb E_{\cl B})
\bb E_{\cl A}\|_{\infty,2}.
\end{align}
\end{pro}

\begin{proof}
Each $x\in K_{\cl A}(e_{\cl B})$ satisfies $0\le x\le 1$, and so the same is
true for elements of $\ovl K^w_{\cl A}(e_{\cl B})$. Moreover, each $x\in
K_{\cl A}(e_{\cl B})$ has unit trace by Lemma \ref{lem2.1} (v). Let $P$ be the
set of finite trace projections in $\langle\cl N, e_{\cl B}\rangle$. Then,
for $x\in \langle\cl N, e_{\cl B}\rangle$, $x\ge 0$,
\begin{equation}\label{eq2.35}
\text{Tr}(x) = \sup\{\text{\rm Tr}(xp)\colon \ p\in P\}.
\end{equation}
If $p\in P$ and $(x_\alpha)$ is a net in $K_{\cl A}(e_{\cl B})$ converging
weakly to $x\in \ovl K^w_{\cl A}(e_{\cl B})$, then
\begin{equation}\label{eq2.36}
\lim_\alpha \text{Tr}(x_\alpha p) = \text{Tr}(xp),
\end{equation}
and it follows from (\ref{eq2.35}) that $\text{Tr}(x) \le 1$. Since $x^2\le x$,
it follows that $x\in L^2(\langle\cl N,e_{\cl B}\rangle,\text{Tr})$. Since
$\text{span}\{P\}$ is  norm dense in 
$L^2(\langle\cl N, e_{\cl B}\rangle,\text{Tr})$, 
we conclude from (\ref{eq2.36}) that $(x_\alpha)$ converges
weakly to $x$ in the Hilbert space. Thus $\ovl K^w_{\cl A}(e_{\cl B})$ is
weakly compact in both $\langle\cl N,e_{\cl B}\rangle$ and $L^2(\langle\cl N,
e_{\cl B}\rangle,\text{Tr})$, and so norm closed in the latter. Thus there
is a unique element $h\in \ovl K^w_{\cl A}(e_{\cl B})$ of minimal
$\|\cdot\|_{\text{Tr},2}$ -- norm.

For each unitary $u\in\cl A$, the map $x\mapsto uxu^*$ is a
$\|\cdot\|_{\text{Tr},2}$ -- norm isometry which leaves $\ovl K^w_{\cl
A}(e_{\cl B})$ invariant. Thus 
\begin{equation}\label{eq2.37}
uhu^* = h,\quad u \text{ unitary in } \cl A,
\end{equation}
by minimality of $h$, so $h\in \cl A' \cap \langle \cl N, e_{\cl B}\rangle$.
This proves (i).

Consider a unitary $u\in\cl A$. Then, by Lemma \ref{lem2.1},
\begin{align}
1 - \text{Tr}(e_{\cl B} ue_{\cl B}u^*) &= 1-\text{Tr}(e_{\cl B}\bb E_{\cl
B}(u)u^*)\nonumber\\
&= 1-\tr(\bb E_{\cl B}(u)u^*)\nonumber\\
&= 1-\tr(\bb E_{\cl B}(u) \bb E_{\cl B}(u)^*)\nonumber\\
&= 1-\|\bb E_{\cl B}(u)\|^2_2\nonumber\\
&= \|(I-\bb E_{\cl B}) (u)\|^2_2\nonumber\\
\label{eq2.38}
&\le \|(I-\bb E_{\cl B}) \bb E_{\cl A}\|^2_{\infty,2}.
\end{align}
This inequality persists when $ue_{\cl B}u^*$ is replaced by elements of
$K_{\cl A}(e_{\cl B})$, so it follows from (\ref{eq2.36}) that
\begin{equation}\label{eq2.39}
1-\text{Tr}(e_{\cl B}h) \le \|(I-\bb E_{\cl B}) \bb E_{\cl A}\|^2_{\infty,2},
\end{equation}
proving (ii).

Since $h\in \cl A'\cap \langle\cl N, e_{\cl B}\rangle$,
\begin{equation}\label{eq2.40}
\text{Tr}(ue_{\cl B}u^*h) = \text{Tr}(e_{\cl B}u^*hu) = \text{Tr}(e_{\cl B}h)
\end{equation}
for all unitaries $u\in\cl A$. Part (iii) follows from this by taking suitable
convex combinations and a weak limit to replace $ue_{\cl B}u^*$ by $h$ on the
left hand side of (\ref{eq2.40}). Finally,  using (\ref{eq2.39}) and
(\ref{eq2.40}),
\begin{align}
\|h-e_{\cl B}\|^2_{\text{Tr},2} &= \text{Tr}(h^2 - 2he_{\cl B} + e_{\cl
B})\nonumber\\ 
&= \text{Tr}(e_{\cl B}-he_{\cl B})\nonumber\\
&= 1 - \text{Tr}(he_{\cl B})\nonumber\\
\label{eq2.41}
&\le \|(I-\bb E_{\cl B})\bb E_{\cl  A}\|^2_{\infty,2}
\end{align}
proving (iv).
\end{proof}

For the last result of this  section, $h$ is the element constructed in the
previous proposition.

\begin{cor}\label{cor2.5}
Let $f$ be the spectral projection of $h$ corresponding to the interval
\newline {\rm [1/2, 1]}. Then $f\in \cl A'\cap\langle\cl N, e_{\cl B}\rangle$,
and \begin{equation}\label{eq2.42}
\|e_{\cl B}-f\|_{\text{\rm Tr},2} \le 2\|(I-\bb E_{\cl B})\bb E_{\cl
A}\|_{\infty,2}. \end{equation}
\end{cor}

\begin{proof}
The first assertion is a consequence of elementary spectral theory. The second
follows from Proposition~\ref{pro2.4} (iv) and Lemma~\ref{lem2.3} (ii).
\end{proof}
\newpage

\section{Estimates in the $\|\cdot\|_2$-norm}\label{sec3}

\indent

This section establishes some more technical results which will be needed
subsequently, the most important of which is Proposition~\ref{pro3.4}. We
maintain the conventions of the previous section:\ $\cl N$ is a separably
acting type $II_1$ factor with masas $\cl A$ and $\cl B$, unless stated to the
contrary.

\begin{lem}\label{lem3.1}
Let $w\in \cl N$ have polar decomposition $w=vk$, where $k = (w^*w)^{1/2}$,
and let $p=v^*v$ and $q=vv^*$ be the initial and final projections of $v$. If
$e\in\cl N$ is a projection satisfying $ew=w$, then
\begin{align}
\label{eq3.1}
\text{\rm (i)}\quad &\|p-k\|_2 \le \|e-w\|_2;\hspace{3.3in}\\
\label{eq3.2}
\text{\rm (ii)}\quad &\|e-q\|_2 \le \|e-w\|_2;\hspace{3.3in}\\
\label{eq3.3}
\text{\rm (iii)}\quad &\|e-v\|_2 \le 2\|e-w\|_2.\hspace{3.3in}
\end{align}
\end{lem}

\begin{proof}
The first inequality is equivalent to
\begin{equation}\label{eq3.4}
\tr(p+k^2 - 2pk) \le \tr(e+w^*w -w-w^*),
\end{equation}
since $ew=w$, and (\ref{eq3.4}) is in turn equivalent to
\begin{equation}\label{eq3.5}
\tr(w+w^*) \le \tr(e-p+2k),
\end{equation}
since $pk=k$ from properties of the polar decomposition. The map $x\mapsto
\tr(k^{1/2}xk^{1/2})$ is a positive linear functional whose norm is $\tr(k)$.
 Thus \begin{equation}\label{eq3.6}
|\tr(w)| = |\tr(w^*)| = |\tr(vk)| = |\tr(k^{1/2}vk^{1/2})| \le \tr(k).
\end{equation}
The range of $e$ contains the range of $w$, so $e\ge q$. Thus
\begin{equation}\label{eq3.7}
\tr(e) \ge \tr(q) = \tr(p),
\end{equation}
and so (\ref{eq3.5}) follows from (\ref{eq3.6}), establishing (i).

The second inequality is equivalent to
\begin{align}
\tr(q-2eq) &\le \tr(w^*w - ew-w^*e)\nonumber\\
\label{eq3.8}
&= \tr(k^2-w-w^*).
\end{align}

Since $eq=q$, this is equivalent to
\begin{equation}\label{eq3.9}
\tr(w+w^*) \le \tr(k^2+q) = \tr(p+k^2).
\end{equation}
From (\ref{eq3.6})
\begin{equation}\label{eq3.10}
\tr(w+w^*)\le 2\ \tr(k) = \tr(p+k^2 - (k-p)^2) \le \tr(p+k^2),
\end{equation}
which establishes (\ref{eq3.9}) and proves (ii).

The last inequality is
\begin{align}
\|e-v\|_2 &\le \|e-vk\|_2 + \|v(k-p)\|_2\nonumber\\
&\le \|e-w\|_2 + \|k-p\|_2\nonumber\\
\label{eq3.11}
&\le 2\|e-w\|_2,
\end{align}
by (i).
\end{proof}

The next result gives some detailed properties of the polar decomposition.

\begin{lem}\label{lem3.2}
Let $\cl A$ be an abelian von Neumann subalgebra of $\cl N$ and let
$\phi\colon \ \cl A\to \cl N$ be a normal $*$-homomorphism. Let $w$ have polar
decomposition
\begin{equation}\label{eq3.12}
w = v(w^*w)^{1/2} = (ww^*)^{1/2}v,
\end{equation}
and let $p = v^*v,\ \ q=vv^*$. If
\begin{equation}\label{eq3.13}
\phi(a) w = wa,\qquad a\in \cl A,
\end{equation}
then
\begin{align}
\label{eq3.14}
\text{\rm (i)}\quad &w^*w\in \cl A' \quad \text{and}\quad ww^*\in \phi(\cl
A)';\hspace{2.5in}\\ 
\label{eq3.15}
\text{\rm (ii)} \quad &\phi(a)v = va\quad \text{and}\quad \phi(a)q = vav^*
\text{ for all}\quad a \in \cl A;\hspace{1in}\\\label{eq3.16}
\text{\rm (iii)}\quad &p\in\cl A'\cap \cl N\quad \text{and}\quad q\in\phi(\cl
A)' \cap \cl N.\hspace{2.2in}
\end{align}
\end{lem}

\begin{proof}
If $a\in\cl A$, then
\begin{align}
w^*wa &= w^*\phi(a)w = (\phi(a^*)w)^*w\nonumber\\
\label{eq3.17}
&= (wa^*)^*w = aw^*w,
\end{align}
and so $w^*w \in\cl A'$. The second statement in (i) has a similar proof.

Let $f$ be the projection onto the closure of the range of $(w^*w)^{1/2}$.
Since $w^* = (w^*w)^{1/2}v^*$, the range of $w^*$ is contained in the range of
$f$, and so $f\ge p$ by Lemma~\ref{lem2.3}(i). For all $x\in\cl N$ and
$a\in\cl A$,
\begin{align}
\phi(a) v(w^*w)^{1/2}x &= \phi(a)wx\nonumber\\
&= wax\nonumber\\
&= v(w^*w)^{1/2}ax\nonumber\\
\label{eq3.18}
&= va(w^*w)^{1/2}x,
\end{align}
since $(w^*w)^{1/2}\in\cl A'$ by (i). Thus
\begin{equation}\label{eq3.19}
\phi(a) vf = vaf,\qquad a\in\cl A,
\end{equation}
which  reduces to
\begin{equation}\label{eq3.20}
\phi(a)v = va,\qquad a\in\cl A,
\end{equation}
since $f\in VN((w^*w)^{1/2}) \subseteq \cl A'$, and
\begin{equation}\label{eq3.21}
v = vp = vpf = vf.
\end{equation}
This proves the first statement in (ii). The second is immediate from
\begin{equation}\label{eq3.22}
\phi(a)q = \phi(a)vv^* = vav^*,\qquad a\in\cl A.
\end{equation}

The proof of the third part is similar to that of the first, and we omit the
details.
\end{proof}

The next result is an approximation to show that certain projections are close
to projections with similar properties, but lying in a given masa.

\begin{pro}\label{pro3.3}
Let $\cl A$ be an abelian von Neumann subalgebra of $\cl N$ and let $\vp>0$.
Let $q\in\cl N$ be a projection satisfying
\begin{itemize}
\item[\rm (i)] $\text{\rm tr}(q)\ge 1-\vp$;
\item[\rm (ii)] $q\in\cl A'\cap \cl N$, and $\cl A q$ is a masa in $q\cl Nq$;
\item[\rm (iii)] if $a\in\cl A$ and $aq=0$, then $a=0$.
\end{itemize}
Let $\cl B$ be a masa containing $\cl A$. Then there exists a projection
$\tilde q \in\cl B$, $\tilde q\le q$, and $\text{\rm tr}(\tilde q)\ge 1-2\vp$.
Moreover, we may choose $\tilde q$ to have the form $pq$ for some projection
$p\in \cl A$.
\end{pro}

\begin{proof}
We will first assume that $\cl A$ contains the identity element of $\cl N$,
and remove this restriction at the end. By (ii), we may define a masa in $q\cl
Bq$ by  $\cl C = \cl Aq$. By (iii), the map
\begin{equation}\label{eq3.23}
aq\mapsto a \mapsto a(1-q),\qquad a\in\cl A,
\end{equation}
is a  well defined $*$-homomorphism $\theta$ of $\cl C$ into $(1-q) \cl
N(1-q)$. Writing elements of $\cl N$ as matrices relative to the decomposition
$1 = q+(1-q)$, we see that 
\begin{equation}\label{eq3.24}
\cl A = \left\{\left(\begin{matrix} c&0\\ 0&\theta(c)\end{matrix}\right)
\colon \ c\in\cl C\right\},
\end{equation}
and $q = \left(\begin{matrix} 1&0\\ 0&0\end{matrix}\right)$.

If 
\begin{equation}\label{eq3.24a}
b = \left(\begin{matrix} b_1&b_2\\ b_3&b_4\end{matrix}\right)
\in \cl B,
\end{equation}
 then $b$ commutes with $\cl A$, so $b_1\in\cl C' = \cl C$, from
(\ref{eq3.24}). Applying this to $b,b^*$ and $bb^*$ gives $b_1,b^*_1$,
$b_2b^*_2\in\cl C$. Thus any partial isometry $v\in q\cl B(1-q)$ (which is
then the (1,2) entry of a matrix in $\cl B$) has the property that
$vv^*\in\cl C$.

Let $\{v_i\}^\infty_{i=1}$ be a set of partial isometries in $q\cl B(1-q)$
which is maximal with respect to having orthogonal initial projections and
orthogonal final projections. Let $q_1$ denote the projection
$\sum\limits^\infty_{i=1} v_iv^*_i\in\cl C$. Normality of the trace implies
that 
\begin{equation}\label{eq3.25}
\tr(q_1) = \sum^\infty_{i=1} \tr(v_iv^*_i) = \tr\left(\sum^\infty_{i=1}
v^*_iv_i\right) \le \vp,
\end{equation}
since each projection $v^*_iv_i$ lies under $1-q$. Let $q_2 = q-q_1$, $q_3 =
\theta(q_1)$, $q_4 = \theta(q_2)$. Then $\{q_i\}^4_{i=1}$ is a set of
orthogonal projections which sum to 1, and we will subsequently write elements
of $\cl N$ as $4\times 4$ matrices relative to the decomposition $1 =
\sum\limits^4_{i=1} q_i$. Define two $*$-homomorphisms $\alpha$ and $\beta$ as
the restrictions of $\theta$ to $\cl Cq_1$ and $\cl Cq_2$ respectively. Then
\begin{equation}\label{eq3.26}
\cl A = \left\{\left(\begin{matrix} f&0&0&0\\ 0&g&0&0\\ 0&0&\alpha(f)&0\\
0&0&0&\beta(g)\end{matrix}\right)\colon \ f\in\cl C q_1, g\in\cl Cq_2\right\}, 
\end{equation}
and $q = \left(\begin{smallmatrix} 1&0&0&0\\ 0&1&0&0\\ 0&0&0&0\\
0&0&0&0\end{smallmatrix}\right)$. Then $q_1+q_3 \in\cl A$ and so commutes with
$\cl B$. This forces each $b\in\cl B$ to have the form
\begin{equation}\label{eq3.25a}
\left(\begin{matrix} *&0&*&0\\ 0&*&0&*\\ *&0&*&0\\
0&*&0&*\end{matrix} \right),
\end{equation}
 so each element of $q\cl B(1-q)$ has the
form 
\begin{equation}\label{eq3.25b}
\left(\begin{matrix} 0&0&*&0\\ 0&0&0&*\\ 0&0&0&0\\
0&0&0&0\end{matrix}\right).
\end{equation}
 Since $q_1v_i = v_i$ for all $i$, it follows
that each $v_i$ has the form 
\begin{equation}\label{3.25c}
\left(\begin{matrix} 0&0&*&0\\ 0&0&0&0\\
0&0&0&0\\ 0&0&0&0\end{matrix}\right),
\end{equation}
 and so
\begin{equation}\label{eq3.27}
v_iv^*_i\le q_1,\quad v^*_iv_i \le q_3,\qquad i\ge 1.
\end{equation}

Suppose that there is an element of $\cl B$ with a non-zero (2,4) or (4,2)
entry. By taking adjoints if necessary, we may assume that the (2,4) entry is
non-zero. Multiplication on the left by $q_2+q_4\in\cl A$ gives an element
\begin{equation}\label{eq3.28}
b = \left(\begin{matrix} 0&0&0&0\\ 0&g&0&x\\ 0&0&0&0\\
0&y&0&z\end{matrix}\right)\in \cl B
\end{equation}
with $g\in\cl Cq_2$ and $x\ne 0$. The condition $\cl B\subseteq \cl A'\cap
\cl N$ leads, as before, to the conclusion that $xx^*$ is a non-zero element
of $\cl Cq_2$ which, being a masa in $q_2\cl Nq_2$, is isomorphic to
$L^\infty[0,1]$. We may then pick a non-negative $h\in\cl Cq_2$ so that $hxx^*$
is a non-zero projection in $\cl Cq_2$. Let $k=h^{1/2}$. Then $v=kx$ is a
non-zero partial isometry. Multiplication of $b$ on the left by
\begin{equation}\label{3.28a}
\left(\begin{matrix} 0&0&0&0\\ 0&k&0&0\\ 0&0&0&0\\
0&0&0&\beta(k)\end{matrix}\right)\in\cl A
\end{equation}
 gives an element
\begin{equation}\label{eq3.28b}\left(\begin{matrix} 0&0&0&0\\ 0&*&0&v\\ 0&0&0&0\\
0&*&0&*\end{matrix}\right)\in\cl B,
\end{equation}
 so 
\begin{equation}\label{eq3.28c}
\left(\begin{matrix}
0&0&0&0\\ 0&0&0&v\\ 0&0&0&0\\ 0&0&0&0\end{matrix}\right)
\in  q\cl B(1-q)
\end{equation}
is a non-zero
partial isometry. The initial and final projections are under
$q_4$ and $q_2$ respectively, and this contradicts the maximality of
$\{v_i\}^\infty_{i=1}$. Thus no (2,4) or (4,2) entry can be non-zero, so each
element of $\cl B$ has the form 
\begin{equation}\label{eq3.28d}
\left(\begin{matrix} *&0&*&0\\ 0&*&0&0\\
*&0&*&0\\ 0&0&0&*\end{matrix}\right).
\end{equation}
 This implies that $q_2\in \cl B'
\cap\cl N = \cl B$. Define $p=q_2+q_4 \in\cl A$, and let $\tilde q = q_2 =
pq$. Then
\begin{equation}\label{eq3.29}
\tr(\tilde q) = \tr(q-q_1)\ge 1-2\vp
\end{equation}
from (i) and (\ref{eq3.25}).

We now turn to the general case. Let $e$ be the unit of $\cl A$ and suppose
that $e\ne 1$. By (ii), $eq=q$ so $e\ge q$. Since $\cl A$ is isomorphic to
$\cl Aq$, by (iii), we may define a normal $*$-isomorphism $\gamma$ from $\cl
A$ to $(1-e)\cl B(1-e)$. We may then define
\begin{equation}\label{eq3.30}
\cl A_1 = \{a+\gamma(a)\colon \ a\in\cl A\},
\end{equation}
which also satisfies the hypotheses. From the first part, there is a
projection $p\in\cl A$ such that
\begin{equation}\label{eq3.31}
\tr((p+\gamma(p))q) \ge 1-2\vp
\end{equation}
and $(p+\gamma(p))q \in\cl B$. Since $\gamma(p)q = 0$, the result follows.
\end{proof}

Recall from \cite{P1} that the normalizing groupoid $\cl G(\cl A)$ of a masa
$\cl A$ in $\cl N$ is the set of partial isometries $v\in\cl N$ such that
$vv^*$, $v^*v\in\cl A$, and $v\cl Av^* = \cl Avv^*$. Such a partial isometry
$v$ implements a spatial $*$-isomorphism between $\cl Av^*v$ and $\cl Avv^*$.
By choosing a normal $*$-isomorphism between the abelian algebras $\cl
A(1-v^*v)$ and $\cl A(1-vv^*)$ (both isomorphic to $L^\infty[0,1]$), we obtain
a $*$-automorphism of $\cl A$ satisfying the hypotheses of Lemma~\ref{lem2.1}
of \cite{JP}. It follows that $v$ has the form $pw^*$, where $p$ is a
projection in $\cl A$ and $w\in N(\cl A)$ (this result is originally in
\cite{Dy}). The next result will allow us to relate $\|\bb E_{\cl A}  - \bb
E_{u\cl Au^*}\|_{\infty,2}$ to the distance from $u$ to $N(\cl A)$.

\begin{pro}\label{pro3.4}
Let $\cl A$ be a masa in $\cl N$, let $u\in\cl N$ be a unitary and let
$\vp_1,\vp_2>0$. Suppose that there exists a partial isometry $v\in\cl N$ such
that $v^*v\in\cl A$, $vv^*\in u\cl A u^*$, $v\cl Av^* = u\cl Au^*vv^*$, and
\begin{gather}
\label{eq3.32}
\|v-\bb E_{u\cl Au^*}(v)\|_2\le \vp_1,\\
\label{eq3.33}
\|v\|^2_2 \ge 1-\vp^2_2.
\end{gather}
Then there exists $\tilde u\in N(\cl A)$ such that
\begin{equation}\label{eq3.34}
\|u-\tilde u\|_2 \le 2(\vp_1+\vp_2).
\end{equation}
\end{pro}

\begin{proof}
Let $v_1$ be the partial isometry $u^*v\in\cl N$. From the hypotheses we see
that $v^*_1v_1$, $v_1v^*_1\in \cl A$ and
\begin{equation}\label{eq3.35}
v_1\cl Av^*_1 = u^*v\cl Av^*u = u^*u\cl A u^*vv^* u = \cl Av_1v^*_1,
\end{equation}
and so $v_1\in \cl G(\cl A)$. It follows from \cite{JP} that $v_1=pw^*$ for
some projection $p\in\cl A$ and unitary $w^*\in N(\cl A)$. Thus
\begin{equation}\label{eq3.36}
vw = up.
\end{equation}
From (\ref{eq3.32}), there exists $a\in\cl A$ such that $\|a\|\le 1$ and $\bb
E_{u\cl Au^*}(v) = uau^*$. Since $\cl A$ is abelian, it is isomorphic to
$C(\Omega)$ for some compact Hausdorff space $\Omega$. Writing $b=|a|$, $0\le
b\le 1$, there exists a unitary $s\in\cl A$ such that $a=bs$.

Now (\ref{eq3.33}) and (\ref{eq3.36}) imply that
\begin{equation}\label{eq3.37}
\|p\|^2_2 =  \|v\|^2_2 \ge 1-\vp^2_2,
\end{equation}
and so
\begin{equation}\label{eq3.38}
\|1-p\|_2 = (1-\|p\|^2_2)^{1/2} \le \vp_2.
\end{equation}
It now follows from (\ref{eq3.36}) that 
\begin{equation}\label{eq3.39}
\|v-uw^*\|_2 = \|vw-u\|_2 = \|up-u\|_2 = \|1-p\|_2 \le \vp_2.
\end{equation}
From (\ref{eq3.32}) and (\ref{eq3.39}) we obtain the estimate
\begin{align}
\|1-bsu^*w\|_2 &= \|uw^* - ubsu^*\|_2\nonumber\\
&\le \|uw^*-v\|_2 + \|v-uau^*\|_2\nonumber\\
\label{eq3.40}
&\le \vp_1+\vp_2.
\end{align}
Let $c = \bb E_{\cl A}(su^*w) \in\cl A$, $\|c\|\le 1$, and apply $\bb E_{\cl
A}$ to (\ref{eq3.40}) to obtain
\begin{equation}\label{eq3.41}
\|1-bc\|_2 \le \vp_1+\vp_2.
\end{equation}
For each $\omega\in\Omega$,
\begin{equation}\label{eq3.42}
|\,1-b(\omega)c(\omega)\,|\ge |\,1-|b(\omega)c(\omega)|\,| \ge 1-b(\omega),
\end{equation}
from which it follows that
\begin{equation}\label{eq3.43}
(1-b)^2 \le (1-bc)(1-bc)^*.
\end{equation}
Apply the trace to (\ref{eq3.43}) and use (\ref{eq3.41}) to reach
\begin{equation}\label{eq3.44}
\|1-b\|_2 \le \vp_1+\vp_2.
\end{equation}
Thus
\begin{equation}\label{eq3.45}
\|a-s\|_2 = \|bs-s\|_2 = \|b-1\|_2 \le \vp_1+\vp_2.
\end{equation}
From (\ref{eq3.32}), (\ref{eq3.44}) and the triangle inequality,
\begin{align}
\|v-usu^*\|_2 &= \|v-ubsu^* + u(b-1)su^*\|_2\nonumber\\
&= \|v-\bb E_{u\cl Au^*}(v) + u(b-1)su^*\|_2\nonumber\\
\label{eq3.46}
&\le 2\vp_1 +\vp_2.
\end{align}
This leads to the estimate
\begin{align}
\|u-ws\|_2 
&= \|su^*w-1\|_2\nonumber\\
&= \|usu^*w-u\|_2\nonumber\\
&= \|usu^*w-up+u(p-1)\|_2\nonumber\\
&\le \|usu^*w-up\|_2 + \vp_2\nonumber\\
&= \|usu^*w-vw\|_2 + \vp_2\nonumber\\
&= \|usu^* - v\|_2 + \vp_2\nonumber\\
\label{eq3.47}
&\le 2(\vp_1+\vp_2),
\end{align}
using (\ref{eq3.38}) and (\ref{eq3.46}). Now define $\tilde u = ws$, which is
in $N(\cl A)$ since $s$ is a unitary in $\cl A$. The last inequality gives
(\ref{eq3.34}).
\end{proof}
\newpage

\section{Homomorphisms on masas}\label{sec4}

\indent

As in the previous sections, $\cl N$ will denote a separably acting type
$II_1$ factor with  masas $\cl A$ and $\cl B$. The $\cl N$-bimodule
$\text{span}\{\cl Ne_{\cl B}\cl N\}$ is weakly dense in $\langle\cl N, e_{\cl
B}\rangle$, \cite{PP}, and there is a well defined $\cl N$-bimodule map
$\Phi\colon \ \text{span}\{\cl Ne_{\cl B}\cl N\}\to \cl N$ which, on 
generators, is given by
\begin{equation}\label{eq4.1}
\Phi(xe_{\cl B}y) = xy,\qquad x,y\in\cl N.
\end{equation}
For subsequent work, the importance of this map is that calculations performed
in $\langle\cl N,e_{\cl B}\rangle$ can be projected down into $\cl N$. The
first lemma collects the properties of $\Phi$ needed later (see \cite{PP,P2}).

\begin{lem}\label{lem4.1}
The map $\Phi\colon \ \text{\rm span}\{\cl Ne_{\cl B}\cl N\}\to \cl N$
satisfies \begin{itemize}
\item[\rm (i)] $\text{\rm tr}(\Phi(x)z) = \text{\rm Tr}(xz)$ for all $x\in
\text{\rm span}\{\cl Ne_{\cl B}\cl N\}$ and for all $z\in\cl N$;
\item[\rm (ii)] $\Phi$ maps $e_{\cl B}\langle\cl N, e_{\cl B}\rangle = e_{\cl
B}\cl N$ into $\cl N$ and
\begin{equation}\label{eq4.3}
\|\Phi(x)\|_2 = \|x\|_{\text{\rm Tr}, 2},\quad x\in e_{\cl B}\langle\cl
N,e_{\cl B} \rangle;
\end{equation}
\item[\rm (iii)] for all $x\in \langle N,e_{\cl B}\rangle$,
\begin{equation}\label{eq4.4}
e_{\cl B}\Phi(e_{\cl B}x) = e_{\cl B}x.
\end{equation}
\end{itemize}
\end{lem}

\begin{proof}
Consider a generator $x = se_{\cl B}t \in \cl Ne_{\cl B}\cl N$. Then, for
$z\in\cl N$,
\begin{equation}\label{eq4.5}
\text{tr}(\Phi(x)z) = \text{tr}(stz) = \text{Tr}(se_{\cl B}tz) = \text{Tr}(xz).
\end{equation}
The first assertion follows by taking linear combinations in (\ref{eq4.5}).

If $x\in e_{\cl B}\cl N$ then write $x = e_{\cl B}y$ for some $y\in\cl N$.
Then $\Phi(x) = y$, so
\begin{equation}\label{eq4.6}
\|x\|^2_{\text{Tr},2} = \text{Tr}(y^*e_{\cl B}y) = \text{tr}(y^*y) =
\|\Phi(x)\|^2_2,
\end{equation}
proving (\ref{eq4.3}).

If $x\in\langle\cl N,e_{\cl B}\rangle$ then $e_{\cl B}x\in e_{\cl B}\cl N$, so
there exists $y\in\cl N$ so that $e_{\cl B}x = e_{\cl B}y$. Then $\Phi(e_{\cl
B}x) =y$, and (\ref{eq4.4}) is immediate.
\end{proof}

The following is the main result of this section.

\begin{thm}\label{thm4.2}
Let $\cl A$ and $\cl B$ be masas in $\cl N$, let $\vp>0$, and let $f$ be a
projection in $\cl A'\cap \langle\cl N,e_{\cl B}\rangle$ with
\begin{equation}\label{eq4.7}
\|e_{\cl B}-f\|_{\text{\rm Tr},2} \le \vp.
\end{equation}
Then there exists a partial isometry $v\in\cl N$ such that
\begin{gather}
\label{eq4.8}
v^*v\in\cl A,\quad vv^*\in\cl B,\quad v\cl Av^*= \cl Bvv^*,\\
\label{eq4.9}
\|v\|^2_2 \ge 1-6\vp^2\\
\intertext{and}
\label{eq4.10}
\|v-\bb E_{\cl B}(v)\|_2 \le (2^{3/2} + 6^{1/2})\vp.
\end{gather}
\end{thm}

\begin{proof}
Various projections will be constructed below. Those associated with $\cl A$
will be denoted by $f,f_0,\ldots, p,p_0,\ldots$ and with $\cl B$ by $e_{\cl
B}, e_0,\ldots, q,q_0,\ldots$~.

By Lemma \ref{lem2.2}, there is a central projection $z\in \langle\cl N,e_{\cl
B}\rangle$ such that $f_1 = zf\in \cl A'\cap \langle\cl N,e_{\cl B}\rangle$
and $e_1 = ze_{\cl B} \le e_{\cl B}$ are equivalent projections in $\langle
\cl N,e_{\cl B}\rangle$. Moreover,
\begin{align}
\|f_1-e_{\cl B}\|^2_{\text{Tr},2} &= \|zf-e_{\cl
B}\|^2_{\text{Tr},2}\nonumber\\
&= \|z(f-e_{\cl B}) - (1-z)e_{\cl B}\|^2_{\text{Tr},2}\nonumber\\
&= \|zf-ze_{\cl B}\|^2_{\text{Tr},2} + \|ze_{\cl B} - e_{\cl
B}\|^2_{\text{Tr},2}\nonumber\\
\label{eq4.11}
&\le 2\vp^2
\end{align}
by (\ref{eq2.12}), and so
\begin{equation}\label{eq4.12}
\|f_1-e_{\cl B}\|_{\text{Tr},2} \le 2^{1/2}\vp.
\end{equation}
Let $v_1\in \langle\cl N,e_{\cl B}\rangle$ be a partial isometry such that 
\begin{equation}\label{eq4.13}
f_1 = v^*_1v_1,\quad e_1=v_1v^*_1
\end{equation}
and
\begin{equation}\label{eq4.14}
\|v_1-e_1\|_{\text{Tr},2}\le 2^{1/2}\|f_1-e_1\|_{\text{Tr},2} \le 2^{1/2}\vp,
\end{equation}
using (\ref{eq2.27}) and (\ref{eq2.12}).

Now define $\theta\colon \ \langle\cl N,e_{\cl B}\rangle \to e_1\langle \cl N,
e_{\cl B}\rangle e_1$ by
\begin{equation}\label{eq4.15}
\theta(x) = v_1xv^*_1,\qquad x\in \langle\cl N,e_{\cl B}\rangle.
\end{equation}
Since $e_{\cl B}e_1 = e_1$ and $e_1v_1 = v_1$, the range of $\theta$ is
contained in $e_{\cl B}\langle\cl N,e_{\cl B}\rangle e_{\cl B} = e_{\cl B}\cl
B$, from Lemma~\ref{lem2.1} (iv). By Lemma~\ref{lem2.1} (iii) there is a well
defined map $\phi\colon \ \langle\cl N,e_{\cl B}\rangle\to \cl B$, given by
\begin{equation}\label{eq4.16}
e_{\cl B} \phi(x) = \theta(x),\qquad x\in\langle \cl N,e_{\cl B}\rangle.
\end{equation}
In particular, $\theta(f_1) = e_1 = e_{\cl B}e_1e_{\cl B}$, and so there is a
projection $e_0 \in\cl B$ so that
\begin{equation}\label{eq4.17}
\theta(f_1) = v_1f_1v^*_1 = e_1 = e_{\cl B}e_0 = e_{\cl B}\phi(f_1).
\end{equation}
Now consider $a_1,a_2\in \cl A$, and note that $f_1\in\cl A'\cap \langle\cl N,
e_{\cl B}\rangle$. Then
\begin{align}
e_{\cl B}\phi(a_1a_2) &= \theta(a_1a_2) = v_1a_1a_2v^*_1\nonumber\\
&= v_1f_1a_1a_2v^*_1 = v_1a_1f_1a_2v^*_1\nonumber\\
&= v_1a_1v^*_1v_1av^*_1 = \theta(a_1)\theta(a_2)\nonumber\\
\label{eq4.17a}
&= e_{\cl B}\phi(a_1) e_{\cl B}\phi(a_2) = e_{\cl B}\phi(a_1)\phi(a_2),
\end{align}
so $\phi(a_1a_2) = \phi(a_1)\phi(a_2)$ by Lemma \ref{lem2.1} (iii). A similar
argument shows that $\phi(a^*) = \phi(a)^*$, and so the restriction of $\phi$
is a normal $*$-homomorphism of $\cl A$ into $\cl B$. Since $v = v_1f_1$ and
$e_{\cl B}v_1=v_1$, we also obtain
\begin{align}
\phi(a)v_1 &= \phi(a) e_{\cl B}v_1 = v_1av^*_1v_1\nonumber\\
\label{eq4.18}
&= v_1af_1 = v_1f_1a = v_1a
\end{align}
for all $a\in\cl A$.

Now let $w = \Phi(v_1)\in \cl N$. Then it follows  from Lemma \ref{lem4.1} that
\begin{equation}\label{eq4.19}
e_{\cl B}w = e_{\cl B}\Phi(v_1) = e_{\cl B}\Phi(e_{\cl B}v_1) = e_{\cl B}v_1 =
v_1.
\end{equation}
Thus, for $a\in\cl A$,
\begin{align}
\phi(a)w &= \phi(a)\Phi(v_1) = \Phi(\phi(a)v_1)\nonumber\\
\label{eq4.20}
&= \Phi(v_1a) = \Phi(v_1)a = wa,
\end{align}
using (\ref{eq4.18}) and the $\cl N$-modularity of $\Phi$. Moreover
\begin{align}
\|w-e_0\|_2 &= \|e_{\cl B}(w-e_0)\|_{\text{Tr},2}\nonumber\\
&= \|v_1-e_1\|_{\text{Tr},2}\nonumber\\
\label{eq4.21}
&\le 2^{1/2}\vp,
\end{align}
using (\ref{eq4.14}). Now let
\begin{equation}\label{eq4.22}
w = v_2(w^*w)^{1/2} = (ww^*)^{1/2} v_2
\end{equation}
be the polar decomposition of $w$ (see Lemma \ref{lem2.3} (i)), and let $p =
v^*_2v_2$, $q = v_2v^*_2$. Equation (\ref{eq4.20}) shows that the hypotheses
of Lemma~\ref{lem3.2} are satisfied and so $p\in \cl A'\cap \cl N$, $q\in
\phi(\cl A)' \cap \cl N$, and
\begin{equation}\label{eq4.23}
\phi(a)v_2 = v_2a,\quad \phi(a)q = v_2av^*_2,\qquad a\in\cl A.
\end{equation}
Since
\begin{align}
e_0w &= e_0\Phi(v_1) = \Phi(e_0v_1)\nonumber\\
&= \Phi(e_0e_{\cl B}v_1) = \Phi(e_1v_1)\nonumber\\
\label{eq4.24}
&=\Phi(v_1) = w,
\end{align}
we may apply Lemma \ref{lem3.1} to $w$ and $e_0$ to obtain
\begin{equation}\label{eq4.25}
\|e_0-q\|_2 \le \|e_0-w\|_2 \le 2^{1/2}\vp,
\end{equation}
where the latter inequality is (\ref{eq4.21}). Putting $a=1$ in (\ref{eq4.23})
gives $e_0q=q$, so $e_0\ge q$. Thus
\begin{equation}\label{eq4.26}
\text{tr}(e_0-q) = \|e_0-q\|^2_2 \le 2\vp^2,
\end{equation}
by (\ref{eq4.25}), and so
\begin{equation}\label{eq4.27}
\text{tr}(q)\ge \text{tr}(e_0) - 2\vp^2.
\end{equation}
Since
\begin{align}
\text{tr}(e_0) &= \text{Tr}(e_{\cl B}e_0) = \text{Tr}(e_1)\nonumber\\
&= \text{Tr}(e_{\cl B}) - \text{Tr}(e_{\cl B}-e_1)\nonumber\\
&= 1-\|e_{\cl B}-e_1\|^2_{\text{Tr},2}\nonumber\\
\label{eq4.28}
&\ge 1-\vp^2,
\end{align}
by (\ref{eq2.12}) with $e_1 = ze_{\cl B}$, these two inequalities combine to
give
\begin{equation}\label{eq4.29}
\text{tr}(q) \ge 1-3\vp^2.
\end{equation}
 
Now $\cl J = \{a\in\cl A\colon \ \phi(a) = 0\}$ is a weakly closed ideal in
$\cl A$, so has the form $\cl J = \cl A(1-f_2)$ for some projection $f_2\in\cl
A$. Then (\ref{eq4.23}) shows that
\begin{equation}\label{eq4.30}
pa = v^*_2v_2a = v^*_2\phi(a)v_2 = 0,\qquad a\in\cl J,
\end{equation}
so $p(1-f_2) = 0$. It follows that $p=pf_2$ and $f_2\ge p$. On the other hand,
(\ref{eq4.22}) shows that $wp=w$ (since $v_2p=v_2$), so multiplication on the
right by $1-p$ in (\ref{eq4.19}) gives
\begin{equation}\label{eq4.31}
v_1(1-p) = 0,\quad e_{\cl B}\phi(1-p) = v_1(1-p)v^*_1 = 0.
\end{equation}
Thus $1-p\in\cl J$, so $1-p\le 1-f_2$. It follows that $p=f_2$ and that
\begin{equation}\label{eq4.32}
\phi(p) = \phi(1) = e_0.
\end{equation}
Now consider an element $a\in\cl A$ such that $\phi(a)q=0$. From
(\ref{eq4.23}),
\begin{equation}\label{eq4.33}
pa = pap = v^*_2v_2av^*_2v_2 = v^*_2\phi(a)qv_2 = 0
\end{equation}
and so $pa=0$. Then $a = (1-p)a \in\cl J$, and thus $\phi(a) = 0$. We conclude
that, for $a\in\cl A$,
\begin{equation}\label{eq4.34}
\phi(a)q = 0 \Rightarrow \phi(a) = 0. 
\end{equation}

Now define $\psi\colon \ p\cl Np\to q\cl Nq$ by
\begin{equation}\label{eq4.35}
\psi(x) = v_2xv^*_2,\qquad x\in p\cl  Np.
\end{equation}
Since $p = v^*_2v_2$ and $q = v_2v^*_2$, it is easy to see that $\psi$ is a
normal $*$-isomorphism. The factor $p\cl Np$ contains a masa $\cl Ap$, and
thus $\psi(\cl Ap)$ is a masa in $q\cl Nq$. From (\ref{eq4.32}) and
(\ref{eq4.23}),
\begin{equation}\label{eq4.36}
\psi(\cl Ap) = v_2\cl Apv^*_2 = \phi(\cl Ap)q = \phi(\cl A) e_0q = \phi(\cl
A)q,
\end{equation}
so $\phi(\cl A)q$ is a masa in $q\cl Nq$. Now consider $a\in\cl A$ and
$b\in\cl B$. Then
\begin{align}
(\phi(a)q)(qbq) &= q\phi(a)bq = qb\phi(a)q\nonumber\\
\label{eq4.37}
&= (qbq)(\phi(a)q),
\end{align}
since $\phi(\cl A) \subseteq \cl B$, and $q\in\phi(\cl A)'$. Thus elements of
$q\cl Bq$ commute with the masa $\phi(\cl A)q$ in $q\cl Nq$, so $q\cl Bq
\subseteq \phi(\cl A)q$. On the other hand,
\begin{equation}\label{eq4.38}
\phi(\cl A)q = q\phi(\cl A)q \subseteq q\cl Bq,
\end{equation}
and this establishes equality.

We have now proved that the hypotheses of Proposition \ref{pro3.3} are
satisfied with $\phi(\cl A)$ replacing $\cl A$ and $3\vp^2$ replacing $\vp$
(see (\ref{eq4.29})). We conclude that there exists a projection $q_1\in
\phi(\cl A)$ such that $q_1q\in\cl B$ and
\begin{equation}\label{eq4.39}
\text{tr}(q_1q) \ge 1-6\vp^2.
\end{equation}
Since $\phi(\cl A) = \phi(\cl Ap)$, there exists a projection $p_1\in\cl A$,
$p_1\le p$, so that $\phi(p_1) = q_1$.

Now define $v=v_2p_1 \in\cl N$. Then
\begin{equation}\label{eq4.40}
v^*v = p_1v^*_2v_2p_1 = p_1pp_1 = p_1\in\cl A,
\end{equation}
and so $v$ is a partial isometry. By (\ref{eq4.23}),
\begin{equation}\label{eq4.41}
vv^* = v_2p_1v^*_2 = \phi(p_1)q = q_1q\in\cl B
\end{equation}
and thus
\begin{equation}\label{eq4.42}
\|v\|^2_2 = \text{tr}(vv^*) = \text{tr}(q_1q) \ge 1-6\vp^2,
\end{equation}
from (\ref{eq4.39}). If $a\in\cl A$ then
\begin{align}
v\cl Av^* &= v_2p_1\cl Av^*_2 = \phi(p_1\cl A)q\nonumber\\
\label{eq4.43}
&= \phi(\cl A)qq_1 = \cl Bvv^*,
\end{align}
by (\ref{eq4.23}) and (\ref{eq4.41}). Finally we estimate $\|v-\bb E_{\cl
B}(v)\|_2$. Since
\begin{align}
\|1-p_1\|^2_2 &=  1-\text{tr}(p_1)\nonumber\\
&\le 1-\text{tr}(v_2p_1v^*_2)\nonumber\\
&= 1-\|v\|^2_2\nonumber\\
\label{eq4.44}
&\le 6\vp^2,
\end{align}
we obtain $\|1-p_1\|_2 \le 6^{1/2}\vp$. Then
\begin{align}
\|v-\bb E_{\cl B}(v)\|_2 &\le \|v-e_0\|_2\nonumber\\
&= \|v_2p_1-e_0\|_2\nonumber\\
&\le \|v_2-e_0\|_2 + \|v_2(1-p_1)\|_2\nonumber\\
&\le 2\|e_0-w\|_2 + 6^{1/2}\vp\nonumber\\
\label{eq4.45}
&\le (2^{3/2} + 6^{1/2})\vp,
\end{align}
where we have used (\ref{eq4.21}) and Lemma \ref{lem3.1} (iii). This completes
the proof.
\end{proof}
\newpage

\section{The main results}\label{sec5}

\indent

In the previous sections, we have completed the necessary preliminary work
and we are now able to prove the main result of the paper,
Theorem~\ref{thm5.3}. This contains two inequalities which we present
separately. We denote by $d(x,S)$ the distance in $\|\cdot\|_2$-norm from an
element $x\in\cl N$ to a subset $S\subseteq \cl N$.

\begin{thm}\label{thm5.1}
Let $\cl A$ be a masa in a separately acting type $II_1$ factor $\cl N$, and
let $u\in\cl N$ be a unitary. Then
\begin{equation}\label{eq5.1}
d(u,N(\cl A)) \le 31 \|(I-\bb E_{u\cl Au^*})\bb E_{\cl A}\|_{\infty,2} \le
31\|\bb E_{\cl A} - \bb E_{u\cl Au^*}\|_{\infty,2}.
\end{equation}
\end{thm}

\begin{proof}
Define $\vp$ to be $\|(I-\bb E_{u\cl Au^*})\bb E_{\cl A}\|_{\infty,2}$. If
$\vp=0$ then $\bb E_{\cl A} = \bb E_{u\cl Au^*}$ and $u\in N(\cl A)$, so there
is nothing to prove. Thus assume $\vp>0$. Let $\cl B = u\cl Au^*$.

By Proposition \ref{pro2.4}, there exists $h\in \cl A'\cap \langle\cl N,e_{\cl
B}\rangle$ satisfying
\begin{equation}\label{eq5.2}
\|h-e_{\cl B}\|_{\text{Tr},2}\le \vp.
\end{equation}
Applying Lemma \ref{lem3.2} (ii), the spectral projection $f$ of $h$
corresponding to the interval [1/2,1] lies in $\cl A'\cap \langle\cl N, e_{\cl
B}\rangle$ and satisfies
\begin{equation}\label{eq5.3}
\|f-e_{\cl B}\|_{\text{Tr},2}\le 2\vp,
\end{equation}
(see Corollary \ref{cor2.5}).
Theorem \ref{thm4.2} (with $\vp$ replaced by $2\vp$) gives the existence of a
partial isometry $v\in\cl N$ satisfying
\begin{gather}
\label{eq5.4}
v^*v\in\cl A,\quad vv^*\in\cl B = u\cl Au^*,\quad v\cl Av^* = \cl Bvv^* = u\cl
Au^*vv^*,\\
\label{eq5.5}
\|v-\bb E_{u\cl Au^*}(v)\|_2 \le (2^{5/2} + (24)^{1/2})\vp,\\
\label{eq5.6}
\|v\|^2_2 \ge 1-24\vp^2.
\end{gather}
We may now apply Proposition \ref{pro3.4}, with $\vp_1 = (2^{5/2} +
(24)^{1/2})\vp$ and $\vp_2 = (24)^{1/2}\vp$, to obtain a normalizing unitary
$\tilde u\in N(\cl A)$ satisfying
\begin{equation}\label{eq5.7}
\|u-\tilde u\|_2 \le 2(2^{5/2} + 2(24)^{1/2})\vp.
\end{equation}
Since $2(2^{5/2} + 2(24)^{1/2}) = 30.90\ldots$, the first inequality follows.
The second is simply
\begin{align}
\|(I-\bb E_{u\cl Au^*})\bb E_{\cl A}\|_{\infty,2} &= \|(\bb E_{\cl A} - \bb
E_{u\cl Au^*})\bb E_{\cl A}\|_{\infty,2}\nonumber\\
\label{eq5.8}
&\le \|\bb E_{\cl A}-\bb E_{u\cl Au^*}\|_{\infty,2},
\end{align}
completing the proof.
\end{proof}	

\begin{lem}\label{thm5.2}
If $\cl A$ is a masa in a type $II_1$ factor $\cl N$ and $u\in\cl N$ is a
unitary, then
\begin{equation}\label{eq5.9}
\|\bb E_{\cl A}-\bb E_{u\cl Au^*}\|_{\infty,2} \le 4d(u, N(\cl A)).
\end{equation}
\end{lem}

\begin{proof}
Let $v\in N(\cl A)$ and define $w$ to be $uv^*$. Then $w\cl Aw^* = u\cl Au^*$,
so it suffices to estimate $\|\bb E_{\cl A}-\bb E_{w\cl Aw^*}\|_{\infty,2}$.
Let $h=1-w$. Then, for $x\in\cl N$, $\|x\|\le 1$,
\begin{align}
\|\bb E_{\cl A}(x) &- \bb E_{w\cl Aw^*}(x)\|_2 = \|\bb E_{\cl A}(x) - w\bb
E_{\cl A}(w^*xw) w^*\|_2\nonumber\\
&= \|w^*\bb E_{\cl A}(x)w-\bb E_{\cl A}(w^*xw)\|_2\nonumber\\
&\le \|w^*\bb E_{\cl A}(x)w - \bb E_{\cl A}(x)\|_2 + \|\bb E_{\cl A}(x) - \bb
E_{\cl A}(w^*xw)\|_2\nonumber\\
&\le \|\bb E_{\cl A}(x)w - w\bb E_{\cl A}(x)\|_2 + \|x-w^*xw\|_2\nonumber\\
&= \|\bb E_{\cl A}(x)h-h\bb E_{\cl A}(x)\|_2 + \|hx-xh\|_2\nonumber\\
&\le 4\|h\|_2\nonumber\\
&= 4\|1-uv^*\|_2\nonumber\\
\label{eq5.10}
&= 4\|v-u\|_2.
\end{align}
Taking the infimum of the right hand side of (\ref{eq5.10}) over all $v\in
N(\cl A)$ gives (\ref{eq5.9}).
\end{proof}

The next theorem summarizes the previous two results.

\begin{thm}\label{thm5.3}
Let $\cl A$ be a masa in a separably acting type $II_1$ factor $\cl N$ and let
$u$ be a unitary in $\cl N$. Then
\begin{equation}\label{eq5.11}
d(u,N(\cl A))/31 \le \|(I-\bb E_{u\cl Au^*})\bb E_{\cl A}\|_{\infty,2} \le
\|\bb E_{\cl A}-\bb E_{u\cl Au^*}\|_{\infty,2} \le 4d(u,N(\cl A)).
\end{equation}
If $\cl A$ is singular, then $\cl A$ is (1/31)-strongly singular.
\end{thm}

\begin{proof}
The inequalities of (\ref{eq5.11}) are proved in Theorem \ref{thm5.1} and
Lemma \ref{thm5.2}. When $\cl A$ is singular, its normalizer is contained in $\cl
A$, so
\begin{equation}\label{eq5.12}
\|u-\bb E_{\cl A}(u)\|_2 \le d(u, N(\cl A))
\end{equation}
holds. Then
\begin{equation}\label{eq5.13}
\|u-\bb E_{\cl A}(u)\|_2 \le 31\|\bb E_{\cl A} - \bb E_{u\cl
Au^*}\|_{\infty,2},
\end{equation}
proving $\alpha$-strong singularity with $\alpha = 1/31$.
\end{proof}

The right hand inequality of (\ref{eq5.11}) is similar to 
\begin{equation}\label{eq5.14}
\|\bb E_{\cl A} - \bb E_{u\cl
Au^*}\|_{\infty,2}\leq 4\|u-\bb E_{\cl A}(u)\|_2,
\end{equation}
which we obtained in \cite[Prop. 2.1]{SS1}, so $u$ being close to $\cl A$ implies
that $\cl A$ and $u{\cl A}u^*$ are also close. We remarked in the introduction 
that there are only two ways in which $\|\bb E_{\cl A} - \bb E_{u\cl
Au^*}\|_{\infty,2}$ can be small, and we now make precise this assertion and justify it.

\begin{thm}\label{thm5.4}
Let $\cl A$ and $\cl B$ be masas in a separably acting type $II_1$ factor $\cl N$,
and let $\delta_1,\,\delta_2,\,\vp > 0$.
\begin{itemize}
\item[\rm (i)] If there are projections $p \in \cl A$, $q \in \cl B$ and a unitary
$u \in \cl N$ satisfying
\begin{equation}\label{eq5.15}
u^*qu=p,\ \ \ u^*q{\cl B}u=p{\cl A},
\end{equation}
\begin{equation}\label{eq5.16}
\|u-{\bb E}_{\cl B}(u)\|_2 \leq \delta_1
\end{equation}
and
\begin{equation}\label{eq5.17}
\text{tr}(p)=\text{tr}(q)\geq 1-{\delta_2}^2,
\end{equation}
then
\begin{equation}\label{eq5.18}
\|\bb E_{\cl A} - \bb E_{\cl
B}\|_{\infty,2}\leq 4\delta_1 + 2\delta_2.
\end{equation}
\item[\rm (ii)] If $\|\bb E_{\cl A} - \bb E_{\cl
B}\|_{\infty,2}\leq \vp$, then there are projections $p \in \cl A$ and $q \in \cl B$,
and a unitary $u \in \cl N$ satisfying
\begin{equation}\label{eq5.19}
u^*qu=p,\ \ u^*q{\cl B}u=p\cl A,
\end{equation}
\begin{equation}\label{eq5.20}
\|u-\bb E_{\cl B}(u)\|_2 \leq (2^{5/2}+2(24)^{1/2})\vp
\end{equation}
and
\begin{equation}\label{eq5.21}
\text{tr}(p)=\text{tr}(q)\geq 1-24\vp^2.
\end{equation}

\end{itemize}
\end{thm}

\begin{proof} 
\noindent (i) Let $\cl C=u^*{\cl B}u$. Then 
\begin{equation}\label{eq5.22}
\|\bb E_{\cl B} - \bb E_{\cl
C}\|_{\infty,2}\leq 4\delta_1,
\end{equation}
 from (\ref{eq5.14}). If $x \in \cl N$ with 
$\|x\|\leq 1$, then 
\begin{align}
\|\bb E_{\cl C}(x)-\bb E_{\cl A}(x)\|_2&\leq\|(1-p)(\bb E_{\cl C}-\bb E_{\cl A})(x)\|_2
+\|p\bb E_{\cl C}(px)-\bb E_{\cl A}(px)\|_2\nonumber\\
\label{eq5.23}&\leq 2\delta_2,
\end{align}
since
\begin{equation}\label{eq5.24}
p{\cl C}=pu^*{\cl B}u=u^*q{\cl B}u=p{\cl A}
\end{equation}
and $p\bb E_{\cl A}(p(\cdot))$ is the projection onto $p{\cl A}$. Then
(\ref{eq5.18}) follows immediately from (\ref{eq5.22}) and (\ref{eq5.23}).

\noindent (ii) As in the proofs of Theorem \ref{thm5.1} and its preceding results
Proposition \ref{pro2.4}, Lemma \ref{lem3.2} and Theorem \ref{thm4.2},
there is a partial isometry $v \in \cl N$ satisfying
\begin{equation}\label{eq5.25}
p=v^*v \in {\cl A},\ \ q=vv^* \in {\cl B},\ \ v^*q{\cl B}v=p{\cl A},
\end{equation}
\begin{equation}\label{eq5.26}
\|v-\bb E_{\cl B}(v)\|_2 \leq (2^{5/2}+(24)^{1/2})\vp
\end{equation}
and
\begin{equation}\label{eq5.27}
\text{tr}(p)=\text{tr}(q)\geq 1 -24\vp^2.
\end{equation}
Let $w$ be a partial isometry which implements the equivalence
\begin{equation}\label{eq5.28}
w^*w=1-p,\ \ \ ww^*=1-q,
\end{equation}
and let $u=v+w$. Then $u$ is a unitary in $\cl N$, since the initial
and final projections of $v$ and $w$ are orthogonal, and
\begin{equation}\label{eq5.29}
u^*q{\cl B}u=v^*q{\cl B}v=p {\cl A}.
\end{equation}
Observe that 
\begin{equation}\label{eq5.30}
\|w-\bb E_{\cl B}(w)\|_2 \leq \|w\|_2 =(\text{tr}(1-p))^{1/2}\leq(24)^{1/2}\vp,
\end{equation}
so that the inequality
\begin{equation}\label{eq5.31}
\|u-\bb E_{\cl B}(u)\|_2\leq (2^{5/2}+2(24)^{1/2})\vp
\end{equation}
follows from (\ref{eq5.26}) and (\ref{eq5.30}).
\end{proof}

\begin{rem}\label{rem5.5}
Recall from \cite{Ch} that Christensen defined ${\cl A}\subset_{\delta} {\cl B}$ to
mean, in our terminology, that $\|(I-\bb E_{\cl B})\bb E_{\cl A}\|_{\infty,2}
\leq \delta$, and then defined the distance between $\cl A$ and $\cl B$ to be
\begin{equation}\label{eq5.32}
\|{\cl A}-{\cl B}\|_2=\text{max}\,\{\|(I-\bb E_{\cl B})\bb E_{\cl A}\|_{\infty,2},\ 
\|(I-\bb E_{\cl A})\bb E_{\cl B}\|_{\infty,2}\}.
\end{equation}
This quantity is clearly bounded by $\|\bb E_{\cl A}-\bb E_{\cl B}\|_{\infty,2}$,
and the reverse inequality 
\begin{equation}\label{eq5.33}
\|\bb E_{\cl A}-\bb E_{\cl B}\|_{\infty,2}\leq 3\|{\cl A}-{\cl B}\|_2
\end{equation}
follows from \cite[Lemma 5.2]{SS1} and the algebraic identity
\begin{equation}\label{eq5.34}
P-Q=P(I-Q)-(I-P)Q,
\end{equation}
valid for all operators $P$ and $Q$. Thus the two notions of distance give equivalent 
metrics on the space of masas in a type $II_1$ factor.$\hfill\square$ 
\end{rem}

We close with a topological result on the space of masas, in the spirit of \cite{Ch,SS1},
which was previously unobtainable.

\begin{cor}\label{cor5.6}
The set of singular masas in a separably acting type $II_1$ factor is closed
in the $\|\cdot\|_{\infty,2}$-metric.
\end{cor}

\begin{proof}
By Theorem \ref{thm5.3}, it suffices to show that those masas, which satisfy
(\ref{eq5.13}) (with any fixed $\alpha>0$ replacing 31) for all unitaries
$u\in \cl N$, form a closed subset. Consider a Cauchy sequence $\{\cl
A_n\}^\infty_{n=1}$ of masas satisfying (\ref{eq5.12}), and fix a unitary
$u\in\cl N$. By \cite{Ch}, the set of masas is closed, so there is a masa $\cl
A$ such that $\lim\limits_{n\to\infty} \|\bb E_{\cl A_n} - \bb E_{\cl
A}\|_{\infty,2} = 0$. Then
\begin{align}
\|u-\bb E_{\cl A}(u)\|_2 &\le \|u-\bb E_{\cl A_n}(u)\|_2 + \|\bb E_{\cl
A_n}(u) - \bb E_{\cl A}(u)\|_2\nonumber\\
&\le \alpha\|\bb E_{u\cl A_nu^*} - \bb E_{\cl A_n}\|_{\infty,2} + \|\bb
E_{\cl A_n} - \bb E_{\cl A}\|_{\infty,2}\nonumber\\
\label{eq5.35}
&\le \alpha\|\bb E_{u\cl A_nu^*} - \bb E_{u\cl Au^*}\|_{\infty,2} + \alpha
\|\bb E_{u\cl Au^*} - \bb E_{\cl A}\|_{\infty,2}
+ \|\bb E_{\cl A_n} - \bb E_{\cl A}\|_{\infty,2},
\end{align}
and the result follows by letting $n\to \infty$.
\end{proof}

\end{document}